\documentclass[12pt]{amsart}
\usepackage{amssymb}
\usepackage{t1enc}
\usepackage[latin2]{inputenc}
\usepackage{verbatim}
\usepackage{amsmath,amsfonts,amssymb,amsthm}
\usepackage[mathcal]{eucal}
\usepackage{enumerate}
\usepackage[centertags]{amsmath}
\usepackage{graphics}

\newtheorem{theorem}{Theorem}
\newtheorem{lemma}{Lemma}

\newtheorem{corollary}{Corollary}

\newtheorem{remark}{Remark}
\numberwithin{equation}{subsection}

\begin{document}
\author{L. E. Persson, G. Tephnadze and P. Wall}
\title[Nörlund means]{Maximal operators of Vilenkin-Nörlund means }
\address{L.-E. Persson, Department of Engineering Sciences and Mathematics,
Lule\aa\ University of Technology, SE-971 87 Lule\aa , Sweden and Narvik
University College, P.O. Box 385, N-8505, Narvik, Norway.}
\email{larserik@ltu.se}
\address{G. Tephnadze, Department of Mathematics, Faculty of Exact and
Natural Sciences, Tbilisi State University, Chavchavadze str. 1, Tbilisi
0128, Georgia}
\email{giorgitephnadze@gmail.com}
\address{P. Wall, Department of Engineering Sciences and Mathematics, Lule%
\aa\ University of Technology, SE-971 87, Lule\aa , Sweden.}
\email{Peter.Wall@ltu.se}
\thanks{The research was supported by Shota Rustaveli National Science
Foundation grant no.13/06 (Geometry of function spaces, interpolation and
embedding theorems}
\date{}
\maketitle

\begin{abstract}
In this paper we prove and discuss some new $\left( H_{p},weak-L_{p}\right) $
type inequalities of maximal operators of Vilenkin-Nörlund means with
monotone coefficients. We also apply these results to prove a.e. convergence
of \ such Vilenkin-Nörlund means. It is also proved that these results are
the best possible in a special sense. As applications, both some well-known
and new results are pointed out.
\end{abstract}



\section{Introduction}

\bigskip\ The definitions and notations used in this introduction can be
found in our next Section. In the one-dimensional case the first result with
respect to the a.e. convergence of Fejér is due to Fine \cite{fine}. Later,
Schipp \cite{Sc} for Walsh series and Pál, Simon \cite{PS} for bounded
Vilenkin series showed that the maximal operator of Fejér means $\sigma
^{\ast }$ \ is of weak type (1,1), from which the a. e. convergence follows
by standard argument \cite{mz}. Fujji \cite{Fu} and Simon \cite{Si2}
verified that $\sigma ^{\ast }$ is bounded from $H_{1}$ to $L_{1}$. Weisz 
\cite{We2} generalized this result and proved boundedness of $\ \sigma
^{\ast }$ \ from the martingale space $H_{p}$ to the space $L_{p},$ for $%
p>1/2$. Simon \cite{Si1} gave a counterexample, which shows that boundedness
does not hold for $0<p<1/2.$ A counterexample for $p=1/2$ was given by
Goginava \cite{Go} (see also \cite{tep2}). In \cite{tep1} it was proved even
stronger result than the maximal operator $\sigma ^{\ast }$ is unbounded. In
fact, it was proved that there exists a martingale $f\in H_{1/2},$ such that
Fejér means of Vilenkin-Fourier series of the martingale $f$ \ are not
uniformly bounded in the space $L_{1/2}.$ Moreover, Weisz \cite{We4} proved
that the following is true:

\textbf{Theorem W1. }The maximal operator of the Fejér means $\sigma ^{\ast
} $ is bounded from the Hardy space $H_{1/2}$ to the space weak-$L_{1/2}$.

Riesz`s logarithmic means with respect to the trigonometric system was
studied by several authors. We mention, for instance, the papers by Sz\`{a}%
sz \cite{Sz} and Yabuta \cite{Ya}. These means with respect to the Walsh and
Vilenkin systems were investigated by Simon \cite{Si1} and Gát \cite{Ga1}.
Blahota and Gát \cite{bg} considered norm summability of Nörlund logarithmic
means and showed that Riesz`s logarithmic means $R_{n}$ have better
approximation properties on some unbounded Vilenkin groups, than the Fejér
means. In \cite{tep5} it was proved that the maximal operator of Riesz`s
means $R^{\ast }$ is not bounded from the Hardy space $H_{1/2}$ to the space 
$weak-L_{1/2}$, but is not bounded from the Hardy space $H_{p}$ to the space 
$L_{p},$ when $0<p\leq 1/2.$ Since the set of Vilenkin polynomials is dense
in $L_{1},$ by well-known density argument due to Marcinkiewicz and Zygmund 
\cite{mz}, we have that $R_{n}f\rightarrow f,$ \ a.e. for \ all \ $f\in
L_{1} $.

Móricz and Siddiqi \cite{Mor} investigated the approximation properties of
some special Nörlund means of Walsh-Fourier series of $L_{p}$ function in
norm. The case when $q_{k}=1/k$ was excluded, since the methods of Móricz
and Siddiqi are not applicable to Nörlund logarithmic means. In \cite{Ga2} Gá%
t and Goginava proved some convergence and divergence properties of the Nö%
rlund logarithmic means of functions in the class of continuous functions
and in the Lebesgue space $L_{1}.$ Among other things, they gave a negative
answer to a question of Móricz and Siddiqi \cite{Mor}. Gát and Goginava \cite%
{Ga3} proved that for each measurable function $\ \phi \left( u\right)
=\circ \left( u\sqrt{\log u}\right) ,$ there exists an integrable function $%
f $ \ such that 
\begin{equation*}
\int_{G_{m}}\phi \left( \left\vert f\left( x\right) \right\vert \right) d\mu
\,\left( x\right) <\infty
\end{equation*}
\ and there exists a set with positive measure, such that the
Walsh-logarithmic means of the function diverges on this set. In \cite{tep4}
it was proved that there exists a martingale $f\in H_{p},$ $(0<p\leq 1),$
such that the maximal operator of Nörlund logarithmic means $L^{\ast }$ is
not bounded in the space $L_{p}.$

In \cite{gog8} Goginava investigated the behaviour of Cesáro means of
Walsh-Fourier series in detail. In the two-dimensional case approximation
properties of Nörlund and Cesáro means was considered by Nagy (see \cite{na}-%
\cite{nagy}). The a.e. convergence of Cesáro means of \ $f\in L_{1}$ was
proved in \cite{gog5}. The maximal operator $\sigma ^{\alpha ,\ast }$ $\
\left( 0<\alpha <1\right) $ of the $\left( C,\alpha \right) $ means of
Walsh-Paley system was investigated by Weisz \cite{we6}. In his paper Weisz
proved that $\sigma ^{\alpha ,\ast }$ is bounded from the martingale space $%
H_{p}$ to the space $L_{p}$ for $p>1/\left( 1+\alpha \right) .$ Goginava 
\cite{gog4} gave a counterexample, which shows that boundedness does not
hold for $0<p\leq 1/\left( 1+\alpha \right) .$ Recently, Weisz and Simon 
\cite{sw} showed the following statement:

\textbf{Theorem SW1. }The maximal operator $\sigma ^{\alpha ,\ast }$ $\left(
0<\alpha <1\right) $ of the $\left( C,\alpha \right) $ means is bounded from
the Hardy space $H_{1/\left( 1+\alpha \right) }$ to the space $%
weak-L_{1/\left( 1+\alpha \right) }$.

In this paper we derive some new $(H_{p},L_{p})$-type inequalities for the
maximal operators of Nörlund means, with monotone coefficients.

The paper is organized as following: In Section 3 we present and discuss the
main results and in Section 4 the proofs can be found. Moreover, in order
not to disturb our discussions in these Sections some preliminaries are
given in Section 2.

\section{Preliminaries}

Denote by $%
\mathbb{N}
_{+}$ the set of the positive integers, $%
\mathbb{N}
:=%
\mathbb{N}
_{+}\cup \{0\}.$ Let $m:=(m_{0,}$ $m_{1},...)$ be a sequence of the positive
integers not less than 2. Denote by 
\begin{equation*}
Z_{m_{k}}:=\{0,1,...,m_{k}-1\}
\end{equation*}
the additive group of integers modulo $m_{k}$.

Define the group $G_{m}$ as the complete direct product of the groups $%
Z_{m_{i}}$ with the product of the discrete topologies of $Z_{m_{j}}`$s.

The direct product $\mu $ of the measures 
\begin{equation*}
\mu _{k}\left( \{j\}\right) :=1/m_{k}\text{ \ \ \ }(j\in Z_{m_{k}})
\end{equation*}%
is the Haar measure on $G_{m_{\text{ }}}$with $\mu \left( G_{m}\right) =1.$

In this paper we discuss bounded Vilenkin groups,\textbf{\ }i.e. the case
when $\sup_{n}m_{n}<\infty .$

The elements of $G_{m}$ are represented by sequences 
\begin{equation*}
x:=\left( x_{0},x_{1},...,x_{j},...\right) ,\ \left( x_{j}\in
Z_{m_{j}}\right) .
\end{equation*}

It is easy to give a base for the neighborhood of $G_{m}:$

\begin{equation*}
I_{0}\left( x\right) :=G_{m},\text{ \ }I_{n}(x):=\{y\in G_{m}\mid
y_{0}=x_{0},...,y_{n-1}=x_{n-1}\},
\end{equation*}%
where $x\in G_{m},$ $n\in 
\mathbb{N}
.$

Denote $I_{n}:=I_{n}\left( 0\right) $ for $n\in 
\mathbb{N}
_{+},$ and $\overline{I_{n}}:=G_{m}$ $\backslash $ $I_{n}$.

\bigskip If we define the so-called generalized number system based on $m$
in the following way : 
\begin{equation*}
M_{0}:=1,\ M_{k+1}:=m_{k}M_{k}\,\,\,\ \ (k\in 
\mathbb{N}
),
\end{equation*}%
then every $n\in 
\mathbb{N}
$ can be uniquely expressed as $n=\sum_{j=0}^{\infty }n_{j}M_{j},$ where $%
n_{j}\in Z_{m_{j}}$ $(j\in 
\mathbb{N}
_{+})$ and only a finite number of $n_{j}`$s differ from zero.

Next, we introduce on $G_{m}$ an orthonormal system which is called the
Vilenkin system. At first, we define the complex-valued function $%
r_{k}\left( x\right) :G_{m}\rightarrow 
\mathbb{C}
,$ the generalized Rademacher functions, by%
\begin{equation*}
r_{k}\left( x\right) :=\exp \left( 2\pi ix_{k}/m_{k}\right) ,\text{ }\left(
i^{2}=-1,x\in G_{m},\text{ }k\in 
\mathbb{N}
\right) .
\end{equation*}

Now, define the Vilenkin system$\,\,\,\psi :=(\psi _{n}:n\in 
\mathbb{N}
)$ on $G_{m}$ as: 
\begin{equation*}
\psi _{n}(x):=\prod\limits_{k=0}^{\infty }r_{k}^{n_{k}}\left( x\right)
,\,\,\ \ \,\left( n\in 
\mathbb{N}
\right) .
\end{equation*}

Specifically, we call this system the Walsh-Paley system when $m\equiv 2.$

The norm (or quasi-norm) of the space $L_{p}(G_{m})$ $\left( 0<p<\infty
\right) $ is defined by 
\begin{equation*}
\left\Vert f\right\Vert _{p}^{p}:=\int_{G_{m}}\left\vert f\right\vert
^{p}d\mu .
\end{equation*}%
\qquad

The space $weak-L_{p}\left( G_{m}\right) $ consists of all measurable
functions $f,$ for which 
\begin{equation*}
\left\Vert f\right\Vert _{weak-L_{p}}^{p}:=\underset{\lambda >0}{\sup }%
\lambda ^{p}\mu \left( f>\lambda \right) <+\infty .
\end{equation*}

The Vilenkin system is orthonormal and complete in $L_{2}\left( G_{m}\right) 
$ (see \cite{Vi}).

Now, we introduce analogues of the usual definitions in Fourier-analysis. If 
$f\in L_{1}\left( G_{m}\right) $ we can define the Fourier coefficients, the
partial sums of the Fourier series, the Dirichlet kernels with respect to
the Vilenkin system in the usual manner:

\begin{equation*}
\widehat{f}\left( n\right) :=\int_{G_{m}}f\overline{\psi }_{n}d\mu \,\ \ \ \
\ \ \,\left( n\in 
\mathbb{N}
\right) ,
\end{equation*}%
\begin{equation*}
S_{n}f:=\sum_{k=0}^{n-1}\widehat{f}\left( k\right) \psi _{k},\text{ \ \ }%
D_{n}:=\sum_{k=0}^{n-1}\psi _{k\text{ }},\text{ \ \ }\left( n\in 
\mathbb{N}
_{+}\right)
\end{equation*}%
respectively.

Recall that 
\begin{equation}
D_{M_{n}}\left( x\right) =\left\{ 
\begin{array}{ll}
M_{n}, & \text{if\thinspace \thinspace \thinspace }x\in I_{n}, \\ 
0, & \text{if}\,\,x\notin I_{n}.%
\end{array}%
\right.  \label{1dn}
\end{equation}

The $\sigma $-algebra generated by the intervals $\left\{ I_{n}\left(
x\right) :x\in G_{m}\right\} $ will be denoted by $\digamma _{n}\left( n\in 
\mathbb{N}
\right) .$ Denote by $f=\left( f^{\left( n\right) },n\in 
\mathbb{N}
\right) $ a martingale with respect to $\digamma _{n}\left( n\in 
\mathbb{N}
\right) .$ (for details see e.g. \cite{We1}).

The maximal function of a martingale $f$ \ is defined by 
\begin{equation*}
f^{\ast }=\sup_{n\in 
\mathbb{N}
}\left\vert f^{(n)}\right\vert .
\end{equation*}

For $0<p<\infty $ \ the Hardy martingale spaces $H_{p}\left( G_{m}\right) $
consist of all martingales for which 
\begin{equation*}
\left\Vert f\right\Vert _{H_{p}}:=\left\Vert f^{\ast }\right\Vert
_{p}<\infty .
\end{equation*}

If $f=\left( f^{\left( n\right) },n\in 
\mathbb{N}
\right) $ is a martingale, then the Vilenkin-Fourier coefficients must be
defined in a slightly different manner: 
\begin{equation*}
\widehat{f}\left( i\right) :=\lim_{k\rightarrow \infty
}\int_{G_{m}}f^{\left( k\right) }\overline{\psi }_{i}d\mu .
\end{equation*}

Let $\{q_{k}:k\geq 0\}$ be a sequence of nonnegative numbers. The $n$-th Nö%
rlund means for a Fourier series of $f$ \ is defined by

\begin{equation}
t_{n}f=\frac{1}{Q_{n}}\overset{n}{\underset{k=1}{\sum }}q_{n-k}S_{k}f,
\label{nor}
\end{equation}%
where $Q_{n}:=\sum_{k=0}^{n-1}q_{k}.$ \ 

We always assume that $q_{0}>0$ and $\ \lim_{n\rightarrow \infty
}Q_{n}=\infty .$\ \ In this case it is well-known that the summability
method generated by $\{q_{k}:k\geq 0\}$ is regular if and only if 
\begin{equation*}
\lim_{n\rightarrow \infty }\frac{q_{n-1}}{Q_{n}}=0.
\end{equation*}

Concerning this fact and related basic results, we refer to \cite{moo}.

If $q_{k}\equiv 1,$ we get the usual Fejér means 
\begin{equation*}
\sigma _{n}f:=\frac{1}{n}\sum_{k=1}^{n}S_{k}f\,.
\end{equation*}%
$\,$

The $\left( C,\alpha \right) $-means of the Vilenkin-Fourier series are
defined by 
\begin{equation*}
\sigma _{n}^{\alpha }f=\frac{1}{A_{n}^{\alpha }}\overset{n}{\underset{k=1}{%
\sum }}A_{n-k}^{\alpha -1}S_{k}f,
\end{equation*}%
where \ 
\begin{equation*}
A_{0}^{\alpha }=0,\text{ \ \ }A_{n}^{\alpha }=\frac{\left( \alpha +1\right)
...\left( \alpha +n\right) }{n!},~~\alpha \neq -1,-2,...
\end{equation*}

The $n$-th Riesz`s logarithmic mean $R_{n}$ and the Nörlund logarithmic mean 
$L_{n}$ are defined by 
\begin{equation*}
R_{n}f:=\frac{1}{l_{n}}\sum_{k=1}^{n-1}\frac{S_{k}f}{k},\text{ \ \ \ }%
L_{n}f:=\frac{1}{l_{n}}\sum_{k=1}^{n-1}\frac{S_{k}f}{n-k},
\end{equation*}%
\ respectively, where $l_{n}:=\sum_{k=1}^{n-1}1/k.$

For the martingale $f$ \ we consider the following maximal operators:%
\begin{equation*}
t^{\ast }f:=\sup_{n\in 
\mathbb{N}
}\left\vert t_{n}f\right\vert ,\text{ \ \ \ }\sigma ^{\ast }f:=\sup_{n\in 
\mathbb{N}
}\left\vert \sigma _{n}f\right\vert ,\text{ \ \ \ \ \ }\sigma ^{\alpha ,\ast
}f:=\sup_{n\in 
\mathbb{N}
}\left\vert \sigma _{n}^{\alpha }f\right\vert ,
\end{equation*}%
\begin{equation*}
R^{\ast }f:=\sup_{n\in 
\mathbb{N}
}\left\vert R_{n}f\right\vert \ \ \ \ \ \ \text{and \ \ \ \ }L^{\ast
}f:=\sup_{n\in 
\mathbb{N}
}\left\vert L_{n}f\right\vert .
\end{equation*}

A bounded measurable function $a$ is a p-atom, if there exists a interval $I$%
, such that%
\begin{equation*}
\int_{I}ad\mu =0,\text{ \ \ }\left\Vert a\right\Vert _{\infty }\leq \mu
\left( I\right) ^{-1/p},\text{ \ \ supp}\left( a\right) \subset I.
\end{equation*}

We also need the following auxiliary results:

\begin{lemma}
\cite{We3} A martingale $f=\left( f^{\left( n\right) },n\in 
\mathbb{N}
\right) $ is in $H_{p}\left( 0<p\leq 1\right) $ if and only if there exists
sequence $\left( a_{k},k\in 
\mathbb{N}
\right) $ of p-atoms and a sequence $\left( \mu _{k},k\in 
\mathbb{N}
\right) ,$ of real numbers, such that, for every $n\in 
\mathbb{N}
,$
\end{lemma}

\begin{equation}
\qquad \sum_{k=0}^{\infty }\mu _{k}S_{M_{n}}a_{k}=f^{\left( n\right) },
\label{1}
\end{equation}%
\begin{equation*}
\sum_{k=0}^{\infty }\left\vert \mu _{k}\right\vert ^{p}<\infty .
\end{equation*}

\textit{Moreover,} 
\begin{equation*}
\left\Vert f\right\Vert _{H_{p}}\backsim \inf \left( \sum_{k=0}^{\infty
}\left\vert \mu _{k}\right\vert ^{p}\right) ^{1/p},
\end{equation*}
\textit{where the infimum is taken over all decomposition of} $f$ \textit{of
the form} (\ref{1}).

\begin{lemma}
\cite{We3} Suppose that an operator $T$ \ is $\sigma $-linear and for some $%
0<p\leq 1$ and
\end{lemma}

\begin{equation*}
\underset{\rho >0}{\sup }\rho ^{p}\mu \left\{ x\in \overline{I}:\left\vert
Ta\right\vert >\rho \right\} \leq c_{p}<\infty ,
\end{equation*}%
\textit{for every} $p$\textit{-atom} $a$\textit{, where} $I$ \ \textit{%
denotes the support of the atom. If} $T$ \textit{is bounded from} $L_{\infty 
\text{ }}$ \textit{to} $L_{\infty },$ \textit{then} 
\begin{equation*}
\left\Vert Tf\right\Vert _{weak-L_{p}}\leq c_{p}\left\Vert f\right\Vert
_{H_{p}}.
\end{equation*}

\textit{and if} $\ 0<p<1,$\textit{\ then} $T$ \ \textit{is of weak
type-(1,1):}%
\begin{equation*}
\left\Vert Tf\right\Vert _{weak-L_{1}}\leq c\left\Vert f\right\Vert _{1}.
\end{equation*}

\section{Main Results}

Our first main result reads:

\begin{theorem}
a) Let $\{q_{k}:k\geq 0\}$ be a sequence of nonnegative numbers, $q_{0}>0$
and $\ $%
\begin{equation*}
\lim_{n\rightarrow \infty }Q_{n}=\infty .
\end{equation*}%
The summability method (\ref{nor}) generated by $\{q_{k}:k\geq 0\}$ is
regular if and only if 
\begin{equation}
\underset{n\rightarrow \infty }{\lim }\frac{q_{n-1}}{Q_{n}}=0.  \label{1a11}
\end{equation}
\end{theorem}

Next, we state our main result concerning the maximal operator of the
summation method (\ref{nor}), which we also show is in a sense sharp.

\begin{theorem}
a) The maximal operator $t^{\ast }$ of \ the summability method (\ref{nor})
with nondecreasing sequence $\{q_{k}:k\geq 0\},$ is bounded from the Hardy
space $H_{1/2}$ to the space $weak-L_{1/2}.$

The statement in a) is sharp in the following sense:

\bigskip b) Let $0<p<1/2$ and\ $\{q_{k}:k\geq 0\}$ is nondecreasing
sequence, satisfying the condition 
\begin{equation}
\frac{q_{0}}{Q_{n}}\geq \frac{c}{n},\text{ \ \ }\left( c>0\right) .
\label{cond1}
\end{equation}%
Then there exists a martingale $f\in H_{p},$ such that 
\begin{equation*}
\underset{n\in 
\mathbb{N}
}{\sup }\left\Vert t_{n}f\right\Vert _{weak-L_{p}}=\infty .
\end{equation*}
\end{theorem}

Our next result shows that the statement in b) above hold also for
nonincreasing sequences and now without any restriction like (\ref{cond1}).

\begin{theorem}
Let $0<p<1/2.$ Then, for all Nörlund means with nonincreasing sequence $%
\left\{ q_{k}:k\geq 0\right\} $, there exists a martingale $f\in H_{p},$
such that 
\begin{equation*}
\underset{n\in 
\mathbb{N}
}{\sup }\left\Vert t_{n}f\right\Vert _{weak-L_{p}}=\infty .
\end{equation*}
\end{theorem}

Up to now we have considered the case $0<p<1/2,$\ but in our final main
result we consider the case when $p=1/\left( 1+\alpha \right) ,$\ $0<\alpha
\leq 1,$\ so that $1/2\leq p<1.$\ Also this result is sharp in two different
important senses.

\begin{theorem}
a) Let $0<\alpha \leq 1.$ Then the maximal operator $t^{\ast }$ of
summability method (\ref{nor}) with non-increasing sequence $\{q_{k}:k\geq
0\},$ satisfying the condition 
\begin{equation}
\frac{n^{\alpha }q_{0}}{Q_{n}}=O\left( 1\right) ,\text{ }\frac{\left\vert
q_{n}-q_{n+1}\right\vert }{n^{\alpha -2}}=O\left( 1\right) ,\text{\ as \ }%
n\rightarrow \infty .  \label{cond3}
\end{equation}%
is bounded from the Hardy space $H_{1/\left( 1+\alpha \right) }$ to the
space $weak-L_{1/\left( 1+\alpha \right) }$.

The parameter $1/\left( 1+\alpha \right) $ in a) is sharp in the following
sense:

\ b) Let $0<p<1/\left( 1+\alpha \right) ,$ where $0<\alpha \leq 1$ and $%
\{q_{k}:k\geq 0\}$ be a non-increasing sequence, satisfying the condition 
\begin{equation}
\frac{q_{0}}{Q_{n}}\geq \frac{c}{n^{\alpha }},\text{ }\left( c>0\right) .
\label{cond4}
\end{equation}%
Then 
\begin{equation*}
\underset{n\in 
\mathbb{N}
}{\sup }\left\Vert t_{n}f\right\Vert _{weak-L_{p}}=\infty .
\end{equation*}

\bigskip Also the condition (\ref{cond3}) is "sharp" in the following sense:

c) Let $\{q_{k}:k\geq 0\}$ be a non-increasing sequence, satisfying the
condition 
\begin{equation}
\overline{\lim_{n\rightarrow \infty }}\frac{q_{0}n^{\alpha }}{Q_{n}}=\infty ,%
\text{ \ \ \ }\left( 0<\alpha \leq 1\right) .  \label{cond2}
\end{equation}%
Then 
\begin{equation*}
\underset{n\in 
\mathbb{N}
}{\sup }\left\Vert t_{n}f\right\Vert _{weak-L_{1/\left( 1+\alpha \right)
}}=\infty .
\end{equation*}
\end{theorem}

A number of special cases of our results are of particular interest and give
both well-known and new information. We just give the following examples of
such Corollaries:

\begin{corollary}
(See \cite{tep4})\textbf{\ }The maximal operator of \ the Fejér means $%
\sigma ^{\ast }$ is bounded from the Hardy space $H_{1/2}$ to the space $%
weak-L_{1/2}$ but is not bounded from $H_{p}$ to the space $weak-L_{p},$
when $0<p<1/2.$
\end{corollary}

\begin{corollary}
(See \cite{tep5}) The maximal operator of \ the Riesz`s means $R^{\ast }$ is
bounded from the Hardy space $H_{1/2}$ to the space $weak-L_{1/2}$ but is
not bounded from $H_{p}$ to the space $weak-L_{p},$ when $0<p<1/2.$
\end{corollary}

\begin{corollary}
The maximal operator of $\ $the $\ \left( C,\alpha \right) $-means $\sigma
^{\alpha ,\ast }$ is bounded from the Hardy space $H_{1/\left( 1+\alpha
\right) }$ to the space $weak-L_{1/\left( 1+\alpha \right) }$ but is not
bounded from $H_{p}$ to the space $weak-L_{p},$ when $0<p<1/\left( 1+\alpha
\right) .$
\end{corollary}

\begin{corollary}
(See \cite{tep4})\textbf{\ }The maximal operator of \ the Nörlund
logarithmic means $L^{\ast }$ is not bounded from the Hardy space $H_{p}$ to
the space $weak-L_{p}$, when $0<p<1.$
\end{corollary}

\begin{corollary}
Let $f\in L_{1}$ and $t_{n}$ be the Nörlund means, with nondecreasing
sequence $\{q_{k}:k\geq 0\}$. Then%
\begin{equation*}
t_{n}f\rightarrow f,\text{ \ \ a.e., \ \ as \ }n\rightarrow \infty .\text{\ }
\end{equation*}
\end{corollary}

\begin{corollary}
Let $f\in L_{1}$ and $t_{n}$ be Nörlund means, with \textit{non-increasing
sequence }$\{q_{k}:k\geq 0\}$ \textit{and} satisfying condition (\ref{cond3}%
). Then 
\begin{equation*}
t_{n}f\rightarrow f,\text{ \ \ a.e., \ \ as \ }n\rightarrow \infty .
\end{equation*}
\end{corollary}

\begin{corollary}
Let $f\in L_{1}$. Then%
\begin{eqnarray*}
\sigma _{n}f &\rightarrow &f,\text{ \ \ \ a.e., \ \ \ as \ }n\rightarrow
\infty ,\text{\ \ \ } \\
R_{n}f &\rightarrow &f,\text{ \ \ \ a.e. \ \ \ \ as \ }n\rightarrow \infty ,
\\
\sigma _{n}^{\alpha }f &\rightarrow &f,\text{ \ \ \ a.e., \ \ as \ }%
n\rightarrow \infty ,\text{ \ \ when \ }0<\alpha <1.
\end{eqnarray*}
\end{corollary}

\begin{remark}
The statements in Corollary 7 are known (see \cite{gog5}, \cite{fine}, \cite%
{sw}, \cite{mz}, \cite{tep9}), but this unified approach to prove them is
new.
\end{remark}

\section{Proofs of the Theorems}

\textbf{Proof of Theorem 1. }The proof is similar as in the case with Walsh
system (see \cite{moo}), so we omit the details.

\textbf{Proof of Theorem 2.} By using Abel transformation we obtain that 
\begin{equation}
Q_{n}:=\overset{n-1}{\underset{j=0}{\sum }}q_{j}=\overset{n}{\underset{j=1}{%
\sum }}q_{n-j}\cdot 1=\overset{n-1}{\underset{j=1}{\sum }}\left(
q_{n-j}-q_{n-j-1}\right) j+q_{0}n,  \label{2b}
\end{equation}%
and

\begin{equation}
t_{n}f=\frac{1}{Q_{n}}\left( \overset{n-1}{\underset{j=1}{\sum }}\left(
q_{n-j}-q_{n-j-1}\right) j\sigma _{j}f+q_{0}n\sigma _{n}f\right) .
\label{2c}
\end{equation}%
Let the sequence $\{q_{k}:k\geq 0\}$ be non-decreasing. By combining (\ref%
{2b}) with (\ref{2c}) and using Abel transformation we get that%
\begin{eqnarray*}
\left\vert t_{n}f\right\vert  &\leq &\left\vert \frac{1}{Q_{n}}\overset{n}{%
\underset{j=1}{\sum }}q_{n-j}S_{j}f\right\vert  \\
&\leq &\frac{1}{Q_{n}}\left( \overset{n-1}{\underset{j=1}{\sum }}\left\vert
q_{n-j}-q_{n-j-1}\right\vert j\left\vert \sigma _{j}f\right\vert
+q_{0}n\left\vert \sigma _{n}f\right\vert \right)  \\
&\leq &\frac{c}{Q_{n}}\left( \overset{n-1}{\underset{j=1}{\sum }}\left(
q_{n-j}-q_{n-j-1}\right) j+q_{0}n\right) \sigma ^{\ast }f\leq c\sigma ^{\ast
}f
\end{eqnarray*}%
so that 
\begin{equation}
t^{\ast }f\leq c\sigma ^{\ast }f.  \label{12aaaa}
\end{equation}

If we apply (\ref{12aaaa}) and Theorem W1 we can conclude that the maximal
operators $t^{\ast }$ are bounded from the Hardy space $H_{1/2}$ to the
space $weak-L_{1/2}.$ It follows that (see Theorem W1) $t^{\ast }$ has weak
type-(1,1) and $t_{n}f\rightarrow f,$ a.e.

In the proof of the second part of Theorem 2 we mainly follow the method of 
\textit{\ }Blahota, Gát and Goginava (see \cite{BGG1}, \cite{BGG2}).

Let $0<p<1/2$ and $\left\{ \alpha _{k}:k\in 
\mathbb{N}
\right\} $ be an increasing sequence of positive integers such that:\qquad 
\begin{equation}
\sum_{k=0}^{\infty }1/\alpha _{k}^{p}<\infty ,  \label{2}
\end{equation}

\begin{equation}
\lambda \sum_{\eta =0}^{k-1}\frac{M_{\alpha _{\eta }}^{1/p}}{\alpha _{\eta }}%
<\frac{M_{\alpha _{k}}^{1/p}}{\alpha _{k}},  \label{3}
\end{equation}

\begin{equation}
\frac{32\lambda M_{\alpha _{k-1}}^{1/p}}{\alpha _{k-1}}<\frac{M_{\alpha
_{k}}^{1/p-2}}{\alpha _{k}},  \label{4}
\end{equation}%
where $\lambda =\sup_{n}m_{n}.$

We note that such an increasing sequence $\left\{ \alpha _{k}:k\in 
\mathbb{N}
\right\} $ which satisfies conditions (\ref{2})-(\ref{4}) can be constructed.

Let \qquad 
\begin{equation}
f^{\left( A\right) }=\sum_{\left\{ k;\text{ }\lambda _{k}<A\right\} }\lambda
_{k}\theta _{k},  \label{55}
\end{equation}%
where

\begin{equation}
\lambda _{k}=\frac{\lambda }{\alpha _{k}}  \label{66}
\end{equation}%
and \ \ \ \ 
\begin{equation}
\theta _{k}=\frac{M_{\alpha _{k}}^{1/p-1}}{\lambda }\left( D_{M_{\alpha
_{k}+1}}-D_{M_{_{\alpha _{k}}}}\right) .  \label{77}
\end{equation}

It is easy to show that the martingale $\,f=\left( f^{\left( 1\right)
},f^{\left( 2\right) }...f^{\left( A\right) }...\right) \in H_{1/2}.$
Indeed, since

\begin{equation*}
S_{M_{A}}\theta _{k}=\left\{ 
\begin{array}{ll}
a_{k}, & \text{if }\ \ \alpha _{k}<A, \\ 
0, & \text{if \ \ }\alpha _{k}\geq A,%
\end{array}%
\right.
\end{equation*}

\begin{equation*}
\text{supp}(\theta _{k})=I_{\alpha _{k}},\text{ \ \ \ \ \ \ \ }%
\int_{I_{\alpha _{k}}}\theta _{k}d\mu =0,\text{ \ \ \ \ }\left\Vert \theta
_{k}\right\Vert _{\infty }\leq M_{\alpha _{k}}^{1/p}=(\text{supp }%
a_{k})^{1/p},
\end{equation*}%
if we apply Lemma 1 and (\ref{2}) we can conclude that $f\in H_{p},$ $\left(
0<p<1/2\right) .$

Moreover, it is easy to show that

\begin{equation}
\widehat{f}(j)=\left\{ 
\begin{array}{ll}
\frac{M_{\alpha _{k}}^{1/p-1}}{\alpha _{k}},\, & \text{if \thinspace
\thinspace }j\in \left\{ M_{\alpha _{k}},...,M_{\alpha _{k}+1}-1\right\} ,%
\text{ }k=0,1,2..., \\ 
0, & \text{if \ \thinspace \thinspace \thinspace }j\notin
\bigcup\limits_{k=1}^{\infty }\left\{ M_{\alpha _{k}},...,M_{\alpha
_{k}+1}-1\right\} .%
\end{array}%
\right.  \label{6}
\end{equation}

We can write

\begin{equation*}
t_{M_{\alpha _{k}}+1}f=\frac{1}{Q_{M_{\alpha _{k}}+1}}\sum_{j=0}^{M_{\alpha
_{k}}}q_{j}S_{j}f+\frac{q_{0}}{Q_{M_{\alpha _{k}}+1}}S_{M_{\alpha
_{k}}+1}f:=I+II.
\end{equation*}

Let $M_{\alpha _{s}}\leq $ $j\leq M_{\alpha _{s}+1},$ where $s=0,...,k-1.$
Moreover, 
\begin{equation*}
\left\vert D_{j}-D_{M_{_{\alpha _{s}}}}\right\vert \leq j-M_{_{\alpha
_{s}}}\leq \lambda M_{_{\alpha _{s}}},\text{ \ \ }\left( s\in 
\mathbb{N}
\right)
\end{equation*}%
so that, according to (\ref{1dn}) and (\ref{6}), we have that%
\begin{eqnarray}
&&\left\vert S_{j}f\right\vert =\left\vert \sum_{v=0}^{M_{\alpha _{s-1}+1}-1}%
\widehat{f}(v)\psi _{v}+\sum_{v=M_{\alpha _{s}}}^{j-1}\widehat{f}(v)\psi
_{v}\right\vert  \label{8} \\
&\leq &\left\vert \sum_{\eta =0}^{s-1}\sum_{v=M_{\alpha _{\eta
}}}^{M_{\alpha _{\eta }+1}-1}\frac{M_{\alpha _{\eta }}^{1/p-1}}{\alpha
_{\eta }}\psi _{v}\right\vert +\frac{M_{\alpha _{s}}^{1/p-1}}{\alpha _{s}}%
\left\vert \left( D_{j}-D_{M_{_{\alpha _{s}}}}\right) \right\vert  \notag \\
&=&\left\vert \sum_{\eta =0}^{s-1}\frac{M_{\alpha _{\eta }}^{1/p-1}}{\alpha
_{\eta }}\left( D_{M_{_{\alpha _{\eta }+1}}}-D_{M_{_{\alpha _{\eta
}}}}\right) \right\vert +\frac{M_{\alpha _{s}}^{1/p-1}}{\alpha _{s}}%
\left\vert \left( D_{j}-D_{M_{_{\alpha _{s}}}}\right) \right\vert  \notag \\
&\leq &\lambda \sum_{\eta =0}^{s-1}\frac{M_{\alpha _{\eta }}^{1/p}}{\alpha
_{\eta }}+\frac{\lambda M_{\alpha _{s}}^{1/p}}{\alpha _{s}}\leq \frac{%
2\lambda M_{\alpha _{s-1}}^{1/p}}{\alpha _{s-1}}+\frac{\lambda M_{\alpha
_{s}}^{1/p}}{\alpha _{s}}\leq \frac{4\lambda M_{\alpha _{k-1}}^{1/p}}{\alpha
_{k-1}}.  \notag
\end{eqnarray}

Let $M_{\alpha _{s-1}+1}+1\leq $ $j\leq M_{\alpha _{s}},$ where $s=1,...,k.$
Analogously to (\ref{8}) we can prove that

\begin{eqnarray*}
&&\left\vert S_{j}f\right\vert =\left\vert \sum_{v=0}^{M_{\alpha _{s-1}+1}-1}%
\widehat{f}(v)\psi _{v}\right\vert =\left\vert \sum_{\eta
=0}^{s-1}\sum_{v=M_{\alpha _{\eta }}}^{M_{\alpha _{\eta }+1}-1}\frac{%
M_{\alpha _{\eta }}^{1/p-1}}{\alpha _{\eta }}\psi _{v}\right\vert \\
&=&\left\vert \sum_{\eta =0}^{s-1}\frac{M_{\alpha _{\eta }}^{1/p-1}}{\alpha
_{\eta }}\left( D_{M_{_{\alpha _{\eta }+1}}}-D_{M_{_{\alpha _{\eta
}}}}\right) \right\vert \leq \frac{2\lambda M_{\alpha _{s-1}}^{1/p}}{\alpha
_{s-1}}\leq \frac{4\lambda M_{\alpha _{k-1}}^{1/p}}{\alpha _{k-1}}.
\end{eqnarray*}

Hence

\begin{equation}
\left\vert I\right\vert \leq \frac{1}{Q_{M_{\alpha _{k}}+1}}%
\sum_{j=0}^{M_{\alpha _{k}}}q_{j}\left\vert S_{j}f\right\vert \leq \frac{%
4\lambda M_{\alpha _{k-1}}^{1/p}}{\alpha _{k-1}}\frac{1}{Q_{M_{\alpha
_{k}}+1}}\sum_{j=0}^{M_{\alpha _{k}}}q_{j}\leq \frac{4\lambda M_{\alpha
_{k-1}}^{1/p}}{\alpha _{k-1}}.  \label{10}
\end{equation}

If we now apply (\ref{6}) and (\ref{8}) we get that 
\begin{eqnarray}
\left\vert II\right\vert  &=&\frac{q_{0}}{Q_{M_{\alpha _{k}}+1}}\left\vert 
\frac{M_{\alpha _{k}}^{1/p-1}}{\alpha _{k}}\psi _{M_{\alpha
_{k}}}+S_{M_{\alpha _{k}}}f\right\vert   \label{100} \\
&=&\frac{q_{0}}{Q_{M_{\alpha _{k}}+1}}\left\vert \frac{M_{\alpha
_{k}}^{1/p-1}}{\alpha _{k}}\psi _{M_{\alpha _{k}}}+S_{M_{\alpha
_{k-1}+1}}f\right\vert   \notag \\
&\geq &\frac{q_{0}}{Q_{M_{\alpha _{k}}+1}}\left( \left\vert \frac{M_{\alpha
_{k}}^{1/p-1}}{\alpha _{k}}\psi _{M_{\alpha _{k}}}\right\vert -\left\vert
S_{M_{\alpha _{k-1}+1}}f\right\vert \right)   \notag \\
&\geq &\frac{q_{0}}{Q_{M_{\alpha _{k}}+1}}\left( \frac{M_{\alpha
_{k}}^{1/p-1}}{\alpha _{k}}-\frac{4\lambda M_{\alpha _{k-1}}^{1/p}}{\alpha
_{k-1}}\right) \geq \frac{q_{0}}{Q_{M_{\alpha _{k}}+1}}\frac{M_{\alpha
_{k}}^{1/p-1}}{4\alpha _{k}}.  \notag
\end{eqnarray}

Without lost the generality we may assume that $c=1$ in (\ref{cond1}). By
combining (\ref{10}) and (\ref{100}) we get 
\begin{eqnarray*}
\left\vert t_{M_{\alpha _{k}}+1}f\right\vert &\geq &\left\vert II\right\vert
-\left\vert I\right\vert \geq \frac{q_{0}}{Q_{M_{\alpha _{k}}+1}}\frac{%
M_{\alpha _{k}}^{1/p-1}}{4\alpha _{k}}-\frac{4\lambda M_{\alpha _{k-1}}^{1/p}%
}{\alpha _{k-1}} \\
&\geq &\frac{M_{\alpha _{k}}^{1/p-2}}{4\alpha _{k}}-\frac{4\lambda M_{\alpha
_{k-1}}^{1/p}}{\alpha _{k-1}}\geq \frac{M_{\alpha _{k}}^{1/p-2}}{8\alpha _{k}%
}.
\end{eqnarray*}

On the other hand%
\begin{equation}
\mu \left\{ x\in G_{m}:\left\vert t_{M_{\alpha _{k}}+1}f\left( x\right)
\right\vert \geq \frac{M_{\alpha _{k}}^{1/p-2}}{8\alpha _{k}}\right\} =\mu
\left( G_{m}\right) =1.  \label{88}
\end{equation}

Let $0<p<1/2.$ Then 
\begin{eqnarray}
&&\frac{M_{\alpha _{k}}^{1/p-2}}{8\alpha _{k}}\cdot \mu \left\{ x\in
G_{m}:\left\vert t_{M_{\alpha _{k}}+1}f\left( x\right) \right\vert \geq 
\frac{M_{\alpha _{k}}^{1/p-2}}{8\alpha _{k}}\right\}  \label{99} \\
&=&\frac{M_{\alpha _{k}}^{1/p-2}}{8\alpha _{k}}\rightarrow \infty ,\text{ \
as }k\rightarrow \infty .  \notag
\end{eqnarray}%
\textbf{\ }

The proof is complete.

\textbf{Proof of Theorem 3. }To prove Theorem 3 we use the martingale (\ref%
{55}), where $\lambda _{k}$ are defined by (\ref{66}), \ for which $\alpha
_{k}$ satisfies conditions (\ref{2})-(\ref{4}) and $\theta _{k}$ are given
by (\ref{77}). It is easy to show that for every non-increasing sequence $%
\{q_{k}:k\geq 0\}$ it automatically holds that 
\begin{equation*}
q_{0}/Q_{M_{\alpha _{k}}+1}\geq 1/\left( M_{\alpha _{k}}+1\right) .
\end{equation*}

By combining (\ref{10}) and (\ref{100}) we see that 
\begin{equation*}
\left\vert t_{M_{\alpha _{k}}+1}f\right\vert \geq \left\vert II\right\vert
-\left\vert I\right\vert \geq \frac{M_{\alpha _{k}}^{1/p-2}}{8\alpha _{k}}.
\end{equation*}

Analogously to\ (\ref{88}) and (\ref{99}) we then get that%
\begin{equation*}
\sup_{k}\left\Vert t_{M_{\alpha _{k}}+1}f\right\Vert _{weak-L_{p}}=\infty .
\end{equation*}

The proof is complete.

\textbf{Proof of Theorem 4. }Let $\{q_{k}:k\geq 0\}$ be a sequence of
non-increasing numbers, which satisfies conditions of theorem 4. In this
case (see \cite{bpt}, \cite{tep9}) it was proved that 
\begin{equation*}
\left\vert F_{n}\right\vert \leq \frac{c\left( \alpha \right) }{n^{\alpha }}%
\overset{\left\vert n\right\vert }{\underset{j=0}{\sum }}M_{j}^{\alpha
}\left\vert K_{M_{j}}\right\vert ,
\end{equation*}%
where 
\begin{equation*}
F_{n}=\frac{1}{Q_{n}}\sum_{k=1}^{n}q_{n-k}D_{k}f
\end{equation*}
is Kernel of Nörlund means$.$

Weisz and Simon \cite{sw} (see also Gát, Goginava \cite{gago}) proved that
the\textbf{\ }maximal operator $\sigma ^{\alpha ,\ast }$ of $\left( C,\alpha
\right) $ $\left( 0<\alpha <1\right) $ means is bounded from the Hardy space 
$H_{1/\left( 1+\alpha \right) }$ to the space $L_{1/\left( 1+\alpha \right)
,\infty }$. Their proof was depend on the following inequality 
\begin{equation*}
\left\vert K_{n}^{\alpha }\right\vert \leq \frac{c\left( \alpha \right) }{%
n^{\alpha }}\sum_{j=0}^{\left\vert n\right\vert }M_{j}^{\alpha }\left\vert
K_{M_{j}}\right\vert ,
\end{equation*}%
where $K_{n}^{\alpha }$ is Kernel of $\left( C,\alpha \right) $ means. Since
our estimation of $F_{n}$ is the same, it is easy to see that the proof of
first part of Theorem 4 will be quiet analogously to the Theorem SW1.

To prove the second part of Theorem 4 we use the martingale (\ref{55}),
where $\lambda _{k}$ is defined by (\ref{66}), \ for which $\alpha _{k}$
satisfies conditions (\ref{2}), (\ref{3}) and 
\begin{equation}
\frac{32\lambda M_{\alpha _{k-1}}^{1/p}}{\alpha _{k-1}}<\frac{M_{\alpha
_{k}}^{1/p-1-\alpha }}{\alpha _{k}},  \label{12}
\end{equation}%
where $\theta _{k}$ is given by (\ref{77}).

We note that such an increasing sequence $\left\{ \alpha _{k}:k\in 
\mathbb{N}
\right\} ,$ which satisfies conditions (\ref{2}), (\ref{3}) and (\ref{12}),
can be constructed.

Let $0\,<p<1/\left( 1+\alpha \right) .$ By combining (\ref{2}) and (\ref{3})
we get that

\begin{equation*}
\left\vert t_{M_{\alpha _{k}}+1}f\right\vert \geq \left\vert II\right\vert
-\left\vert I\right\vert =\frac{M_{\alpha _{k}}^{1/p-1}}{4\alpha _{k}}\frac{%
q_{0}}{Q_{M_{\alpha _{k}}+1}}-\frac{4\lambda M_{\alpha _{k-1}}^{1/p}}{\alpha
_{k-1}}.
\end{equation*}

Without lost the generality we may assume that $c=1$ in (\ref{cond4}). Since 
$1/p-1-\alpha >0$ by using (\ref{12}) we find that%
\begin{equation*}
\left\vert t_{M_{\alpha _{k}}+1}f\right\vert \geq \frac{M_{\alpha
_{k}}^{1/p-1-\alpha }}{8\alpha _{k}}
\end{equation*}%
and 
\begin{eqnarray*}
&&\frac{M_{\alpha _{k}}^{1/p-1-\alpha }}{8\alpha _{k}}\cdot \mu \left\{ x\in
G_{m}:\left\vert t_{M_{\alpha _{k}}+1}f\right\vert \geq \frac{M_{\alpha
_{k}}^{1/p-1-\alpha }}{8\alpha _{k}}\right\} \\
&=&\frac{M_{\alpha _{k}}^{1/p-1-\alpha }}{8\alpha _{k}}\rightarrow \infty ,%
\text{ \ as }k\rightarrow \infty .
\end{eqnarray*}

\bigskip Finally, we prove the third part of Theorem 4. Under condition (\ref%
{cond2}) there exists an increasing sequence $\left\{ \alpha _{k}:k\in 
\mathbb{N}
\right\} ,$ which satisfies the conditions 
\begin{equation}
\sum_{k=0}^{\infty }\left( \frac{Q_{M_{\alpha _{k}}+1}}{q_{0}M_{\alpha
_{k}}^{\alpha }}\right) ^{1/2\left( 1+\alpha \right) },  \label{33}
\end{equation}

\begin{equation}
\lambda \sum_{\eta =0}^{k-1}\frac{Q_{M_{\alpha _{\eta }}+1}M_{\alpha _{\eta
}}^{\alpha /2+1}}{q_{0}^{1/2}}\leq \frac{Q_{M_{\alpha _{k}}+1}M_{\alpha
_{k}}^{\alpha /2+1}}{q_{0}^{1/2}},  \label{34}
\end{equation}

\begin{equation}
\frac{32\lambda Q_{M_{\alpha _{k-1}}+1}M_{\alpha _{k-1}}^{\alpha /2+1}}{%
q_{0}^{1/2}}<\left( \frac{q_{0}M_{\alpha _{k}}^{\alpha }}{Q_{M_{\alpha
_{k}}+1}}\right) ^{1/2},  \label{35}
\end{equation}%
where $\lambda =\sup_{n}m_{n}.$

We note that such an increasing sequence $\left\{ \alpha _{k}:k\in 
\mathbb{N}
\right\} ,$ which satisfies conditions (\ref{33})-(\ref{35}), can be
constructed.

Let \qquad 
\begin{equation*}
f^{\left( A\right) }=\sum_{\left\{ k;\text{ }\lambda _{k}<A\right\} }\lambda
_{k}\theta _{k},
\end{equation*}%
where

\begin{equation*}
\lambda _{k}=\lambda \cdot \left( \frac{Q_{M_{\alpha _{k}}+1}}{%
q_{0}M_{\alpha _{k}}^{\alpha }}\right) ^{1/2}
\end{equation*}%
and $\theta _{k}$ are given by (\ref{77}) for $p=1/\left( 1+\alpha \right) .$
If we apply Lemma 1 and (\ref{33}) analogously to the proof of the second
part of Theorem 2 we can conclude that $\,f=\left( f^{\left( 1\right)
},f^{\left( 2\right) }...f^{\left( A\right) }...\right) \in H_{1/\left(
1+\alpha \right) }.$

It is easy to show that

\begin{equation}
\widehat{f}(j)=\left\{ 
\begin{array}{ll}
\left( \frac{Q_{M_{\alpha _{k}}+1}M_{\alpha _{k}}^{\alpha }}{q_{0}}\right)
^{1/2},\, & \text{if \thinspace \thinspace }j\in \left\{ M_{\alpha
_{k}},...,M_{\alpha _{k}+1}-1\right\} ,\text{ }k=0,1,2..., \\ 
0, & \text{if \ \thinspace \thinspace \thinspace }j\notin
\bigcup\limits_{k=1}^{\infty }\left\{ M_{\alpha _{k}},...,M_{\alpha
_{k}+1}-1\right\} ,%
\end{array}%
\right.  \label{26}
\end{equation}%
and

\begin{equation*}
t_{M_{\alpha _{k}}+1}f=\frac{1}{Q_{M_{\alpha _{k}}+1}}\sum_{j=0}^{M_{\alpha
_{k}}}q_{j}S_{j}f+\frac{q_{0}}{Q_{M_{\alpha _{k}}+1}}S_{M_{\alpha
_{k}}+1}f:=III+IV.
\end{equation*}

Let $M_{\alpha _{s}}\leq $ $j\leq M_{\alpha _{s}+1},$ where $s=0,...,k-1.$
Analogously to (\ref{8}) from (\ref{26}) we have that%
\begin{eqnarray}
\left\vert S_{j}f\right\vert &\leq &\left\vert \sum_{\eta =0}^{s-1}\left( 
\frac{Q_{M_{\alpha _{\eta }}+1}M_{\alpha _{\eta }}^{\alpha }}{q_{0}}\right)
^{1/2}\left( D_{M_{_{\alpha _{\eta }+1}}}-D_{M_{_{\alpha _{\eta }}}}\right)
\right\vert  \label{89} \\
&&+\left( \frac{Q_{M_{\alpha _{s}}+1}M_{\alpha _{s}}^{\alpha }}{q_{0}}%
\right) ^{1/2}\left\vert \left( D_{j}-D_{M_{_{\alpha _{s}}}}\right)
\right\vert  \notag \\
&\leq &\frac{4\lambda Q_{M_{\alpha _{k-1}}+1}^{1/2}M_{\alpha _{k-1}}^{\alpha
/2+1}}{q_{0}^{1/2}}.  \notag
\end{eqnarray}

Let $M_{\alpha _{s-1}+1}+1\leq $ $j\leq M_{\alpha _{s}},$ where $\
s=1,...,k. $ Then

\begin{eqnarray*}
\left\vert S_{j}f\right\vert &=&\left\vert \sum_{\eta =0}^{s-1}\left( \frac{%
Q_{M_{\alpha _{\eta }}+1}M_{\alpha _{\eta }}^{\alpha }}{q_{0}}\right)
^{1/2}\left( D_{M_{_{\alpha _{\eta }+1}}}-D_{M_{_{\alpha _{\eta }}}}\right)
\right\vert \\
&\leq &\lambda \sum_{\eta =0}^{s-1}\frac{Q_{M_{\alpha _{\eta }}+1}M_{\alpha
_{\eta }}^{\alpha /2+1}}{q_{0}^{1/2}}\leq \frac{2\lambda Q_{M_{\alpha
_{s-1}}+1}M_{\alpha _{s-1}}^{\alpha /2+1}}{q_{0}^{1/2}}\leq \frac{4\lambda
Q_{M_{\alpha _{k-1}}+1}M_{\alpha _{k-1}}^{\alpha /2+1}}{q_{0}^{1/2}},
\end{eqnarray*}%
and

\begin{equation}
\left\vert III\right\vert \leq \frac{1}{Q_{M_{\alpha _{k}}+1}}%
\sum_{j=0}^{M_{\alpha _{k}}}q_{j}\left\vert S_{j}f\right\vert \leq \frac{%
4\lambda Q_{M_{\alpha _{k-1}}+1}M_{\alpha _{k-1}}^{\alpha /2+1}}{q_{0}^{1/2}}%
\frac{1}{Q_{M_{\alpha _{k}}+1}}\sum_{j=0}^{M_{\alpha _{k}}}q_{j}\leq \frac{%
4\lambda Q_{M_{\alpha _{k-1}}+1}M_{\alpha _{k-1}}^{\alpha /2+1}}{q_{0}^{1/2}}%
.  \label{101}
\end{equation}

If we apply (\ref{26}) and (\ref{89}) we get that%
\begin{eqnarray}
\left\vert IV\right\vert &\geq &\left( \frac{Q_{M_{\alpha _{k}}+1}M_{\alpha
_{k}}^{\alpha }}{q_{0}}\right) ^{1/2}\frac{q_{0}\left\vert D_{M_{\alpha
_{k}}+1}-D_{M_{\alpha _{k}}}\right\vert }{Q_{M_{\alpha _{k}}+1}}-\frac{q_{0}%
}{Q_{M_{\alpha _{k}}+1}}\left\vert S_{M_{\alpha _{k-1}+1}}f\right\vert
\label{102} \\
&\geq &\left( \frac{q_{0}M_{\alpha _{k}}^{\alpha }}{Q_{M_{\alpha _{k}}+1}}%
\right) ^{1/2}-\frac{4\lambda Q_{M_{\alpha _{k-1}}+1}M_{\alpha
_{k-1}}^{\alpha /2+1}}{q_{0}^{1/2}}\geq \frac{1}{4}\left( \frac{%
q_{0}M_{\alpha _{k}}^{\alpha }}{Q_{M_{\alpha _{k}}+1}}\right) ^{1/2}.  \notag
\end{eqnarray}

By combining (\ref{35}), (\ref{101}) and (\ref{102}) we get that 
\begin{eqnarray*}
\left\vert t_{M_{\alpha _{k}}+1}f\right\vert &\geq &\left\vert IV\right\vert
-\left\vert III\right\vert \geq \frac{1}{4}\left( \frac{q_{0}M_{\alpha
_{k}}^{\alpha }}{Q_{M_{\alpha _{k}}+1}}\right) ^{1/2}-\frac{4\lambda
Q_{M_{\alpha _{k-1}}+1}M_{\alpha _{k-1}}^{\alpha /2+1}}{q_{0}^{1/2}} \\
&\geq &\frac{1}{8}\left( \frac{q_{0}M_{\alpha _{k}}^{\alpha }}{Q_{M_{\alpha
_{k}}+1}}\right) ^{1/2}.
\end{eqnarray*}

Hence, it yields that%
\begin{eqnarray*}
&&\frac{1}{8}\left( \frac{q_{0}M_{\alpha _{k}}^{\alpha }}{Q_{M_{\alpha
_{k}}+1}}\right) ^{1/2\left( 1+\alpha \right) }\mu \left\{ x\in
G_{m}:\left\vert t_{M_{\alpha _{k}}+1}f\left( x\right) \right\vert \geq 
\frac{1}{8}\left( \frac{q_{0}M_{\alpha _{k}}^{\alpha }}{Q_{M_{\alpha _{k}}+1}%
}\right) ^{1/2}\right\} \\
&=&\frac{1}{8}\left( \frac{q_{0}M_{\alpha _{k}}^{\alpha }}{Q_{M_{\alpha
_{k}}+1}}\right) ^{1/2}\mu \left( G_{m}\right) =\frac{1}{8}\left( \frac{%
q_{0}M_{\alpha _{k}}^{\alpha }}{Q_{M_{\alpha _{k}}+1}}\right)
^{1/2}\rightarrow \infty ,\text{ \ as \ }k\rightarrow \infty .
\end{eqnarray*}%
\textbf{\ }

The proof is complete.

\textbf{A final remark: }\ Several of the operators considered in this
paper, e.g. those described by the Nörlund means, Riesz`s logarithmic means
and Nörlund logarithmic means are called Hardy type operators in the
literature. The mapping properties of such operators, especially between
weighted Lebegue spaces, is much studied in the literature, see e.g. the
books \cite{kp} and \cite{mkp} and the references there. Such complimentary
information can be of interest for further studies of the problems
considered in this paper.

\textbf{Acknowledgment: }We thank the careful referees for some good
suggestions, which improved the final version of this paper.

\end{document}